\mag=1200
\hsize 16 truecm \hoffset 0.46 truecm
\vsize 23.6truecm \voffset -0.04 truecm
\baselineskip=16pt plus 1pt
\parskip=5pt plus 1pt
\overfullrule = 0pt

\def\Diag{\def\normalbaselines{\baselineskip25pt\lineskip3pt
        \lineskiplimit3pt}%
        \matrix}
\def\hfleche#1#2{\smash{\mathop{\vbox{\hbox to 8.5mm{{#1}}}}\limits^{#2}}}
\def\vfleche#1#2{\bigg #1 \rlap{$\vcenter{\hbox{$\scriptstyle#2$}}$}}

\def\pe#1{{\bar #1}}
\def\p{\pe{p}}
\def\q{\pe{q}}
\def\t{\pe{t}}
\def\0{\pe{0}}
\def\RR{\bf R}

\def\Rm{{\bf R}^*_+}
\def\jbar{\bar {\hbox{\j }}}

\font\sevensl=cmsl10 scaled 800
\font\moyen=cmr10 scaled 800

\centerline{\bf DIFFERENTIAL FORMS ON SINGULAR VARIETIES}
\centerline{\bf  AND CYCLIC HOMOLOGY.}
\centerline{\sl Preliminary version. 20 Nov. 1996}
\centerline{\sl by J. P. Brasselet and A. Legrand}
\medskip

\noindent {\bf Abstract.} A classical result of A. Connes asserts
that the Frechet algebra of smooth functions on a smooth compact manifold
$X$ provides, by a purely algebraic procedure, the de Rham cohomology of
$X$. Namely the procedure uses Hochschild and cyclic homology of this
algebra.

In the situation of a Thom-Mather stratified variety, we construct a
Frechet algebra of functions on the regular part and a module of poles
along the singular part. We associate to these objects a complex  of
differential forms and an Hochschild complex, on the regular part,
both with poles along the singular part. The de Rham cohomology of the
first complex and the cylic homology of the second one are related to
the intersection homology of the variety, the corresponding perversity
is determined by the orders of poles.

The detailed proofs of the results will appear in a forthcoming
publication. The first author thanks the organizers of the Terry Wall's
60th birthday meeting in Liverpool for the invitation to give a lecture
during the conference.

\bigskip

\noindent{\bf 1. Introduction.}

The aim of this paper is to extend to singular varieties the de Rham and
Connes theorems, using the intersection homology, due to Goresky and
MacPherson which is more adapted to singularities, instead of the
classical homology. Let $X$ be a  $C^\infty$ compact manifold,
$\Omega^*(X)$ the de Rham complex and $H^*_{\rm dR}(X)$ its cohomology.
A. Connes shows that the classical differential de Rham construction
$$X     \Longrightarrow \Omega^*(X)\Longrightarrow H^*_{\rm dR}(X)$$
can be provided by a purely algebraic process
$$C^\infty(X) \Longrightarrow
HH^*(C^\infty(X))\Longrightarrow PHC^*(C^\infty(X))$$
where $HH^*(C^\infty(X))$ is the Hochschild homology of the
Frechet algebra of  ${\cal C}^\infty$-functions on $X$ and
$PHC^*(C^\infty(X))$ is the periodic cyclic homology [Co].

Now let $X$ be a Thom-Mather stratified variety, with singular part
$\Sigma \subset X $, so the regular stratum  $X-\Sigma$ is an open
manifold. We want to construct an algebra  $IC^\infty(X)$ of
differentiable functions on $X-\Sigma$ such that
$$C^\infty_c(X-\Sigma )\subset  IC^\infty(X)\subset C^\infty(X-\Sigma)$$
and which contains enough geometric informations on the compactification
of $X-\Sigma$ by the singular variety $\Sigma$ to fullfill the Connes
program in this singular framework. Evidently the two extrem algebras do
not fit in, we need to add asymptotic control to keep informations ``near
$\Sigma$". We shall proceed of the following way~: For each stratum $S_i$
in $\Sigma$, there is a tubular neighborhood $T_i$ of $S_i$ and a
distance function $r_i$ to $S_i$, defined in $T_i$. Firstly we define
the (Frechet) algebra $A$ of differentiable functions in the variables
$r_i$, which are indefinitly logarithmicaly controlled near $\Sigma$.
This will be our asymptotic control ``reference algebra". To be
controlled a space should be an $A$-space, so  $A$ will be also our
``basic ring" for the homological constructions. Then we define an
$A$-algebra of controlled differentiable functions on $X-\Sigma $ and an
$A$-differentiable module of poles along $\Sigma $. These objects allow
us to construct two complexes said $\bar\beta $-controlled.\par

The first one is a complex of differential forms with coefficients in
the module of poles and whose cohomology is the Goresky-MacPherson
intersection cohomology (Theorem 1). The perversity is related to the
orders of poles. The second one is a mixed complex copied from the
Connes's algebraic procedure. Its Hochschild homology is identified with
the first complex and its periodic cyclic homology is the intersection
cohomology of $X$. This last result is explicited in the case of the
cone (Theorem 3). \par

Let us mention the basic difficulties appearing in Hochschild and cyclic
homology in the singular situation, using the intersection homology and
a ``control near $\Sigma $". At the Hochschild homology level~:

(1) Whatever the way to introduce a control (by use of poles or
$L^p$-forms for instance), it is not compatible with the product.

\noindent At the cyclic homology level~:

(2) It is clear that any type of control does not agree with the de
Rham differential (this is the case for the intersection complex  or the
$L^p$-forms complex).

(3) We need to use an Hochschild complex with coefficients but the
cyclic structure does not exist in this case.

To solve these problems we use (cf. \S 4)~:

(1) a specific $A$-equivariant Hochschild theory with
coefficients in poles associated to a modified differentiable
Frechet structure,

(2) a suitable unitarization (theorem 2),

(3) the more general setting  given by mixed complexes.

Although we often speak of $A$-objects, our equivariant  construction
 does not agree with equivariant theories developped in [Bry], [BG].
We do not know if they  can be adapted to this singular situation and
roughly speaking we substitute associated equivariant conditions by
adapted semi-norms. We will precise this point in the paragraph 3 (see
the ``the cone situation").\par  \bigskip

\noindent{\bf 2. Definition of intersection homology.}
\medskip

We will consider  a singular variety
$X$ endowed with a Thom-Mather $C^\infty$ stratification i.e. a
filtration of $X$ by closed subsets
$$X=X_n \supset X_{n-1}\supset X_{n-2} \supset \cdots \supset X_1 \supset
X_0 \supset  X_{-1} = \emptyset  \eqno(*)$$
where $\Sigma =X_{n-1} $ is the singular part. Each stratum $S_i =
X_i-X_{i-1}$ is either an emptyset or an  $i$-dimensional
$C^\infty$-manifold and there are \hfill\break
\hglue 1truecm an open neighborhood $T_i$ of $S_i$ in $X$,\hfill\break
\hglue 1truecm a continuous retraction $\pi_i$ of $T_i$ on
$S_i$,\hfill\break
\hglue 1truecm a continuous function $\rho_i : T_i \to [0,1[$,\hfill
\break such that $S_i = \{ x\in T_i \vert \rho_i(x) =0 \}$ and the
$(T_i,\pi_i,\rho_i)$ satisfy the axioms of Mather [Ma].

These data imply the following local triviality condition~:
$$\forall \ x\in
S_{n-j},\quad \exists \  U_x\subset X \hbox{ and an homeomorphism }
\psi_x:U_x\to B^{n-j}\times cL_x$$ where $B^{n-j}$ is the standard
open $(n-j)$-dimensional ball and $cL_x$ is the cone over the link $L_x$.
The link is assumed to be stratified and independant
of the point $x\in S_{n-j}$ and $\psi_x$ preserves the stratifications of
$U_x$ (induced by the one of $X$) and the one of the product $B^{n-j}
\times cL_x$. The parameter of the cone corresponds to the Mather
distance function $\rho_{n-j}$. For a complete definition see for
instance [GM2].

Now let us recall the definition of intersection homology due to
M. Goresky and R. MacPherson [GM1]. Given a stratified singular variety,
the idea of intersection homology is to consider chains and cycles whose
intersections with the strata are ``not too big". The allowed  chains
and cycles meet the strata with a controlled and fixed defect of
transversality. This defect is an integer function, called a {\sl
perversity}, increasing with the codimension $j$ of the strata, and
denoted $\bar p = (p_0,p_1,p_2,\ldots , p_j,\ldots , p_n)$. It
satisfies~: $$p_0 = p_1 = p_2 =0 \qquad
\hbox{ and } \qquad p_{j}\le p_{j+1} \le p_{j} +1 \quad \hbox{
for } j\ge 2$$
Let $C_i(X)$ be any ``classical" chain complex on $X$ with integer
coefficients, we can define the complex~:
$$\eqalign{IC^{\bar p}_i(X) = \{ \xi \in C_i(X) \ : \
&\dim (\vert \xi \vert \cap X_{n-j})\le  i-j +p_{j} \quad\hbox{ and }\cr
&\dim (\vert \partial \xi \vert \cap X_{n-j})\le  i-1-j +p_{j} \}\ ,\cr}
$$
the {\sl intersection homology groups} $IH^{\p}_i(X)$ are homology groups
of this  complex.

For the zero perversity (all $p_j$ are zero), allowed chains and cycles
are transverse to all strata. The {\sl total} perversity $\t$ is the one
such that, for all $j\ge 2$, $t_j = j-2$. \par
\medskip

We shall be mainly interested by the axiomatic definition of
 intersection homology. Namely, if a complex of sheaves on $X$
satisfies the so-called {\sl perverse sheaves} axioms [GM2], then the
hypercohomology of ($X$ with value in) this perverse sheaf is the
intersection homology of $X$. The main axioms of perverse sheaves are
issued from the following local computation property (cf [GM2]).

Let $cL$ be the open cone over an $(n-1)$-dimensional manifold $L$, then
the perversity depends only on $p_n$ and we have~:
$$IH^{\p}_i({c}(L))\cong \cases{
H_i(L) & $i<n-p_n-1$\cr
{}&{}\cr
0 & $i\ge n-p_n-1$\cr}$$

The intersection homology is the good theory for extending many of
classical results from manifolds to singular varieties. \par

The most important is  Poincar\'e  duality (which motives the theory).
The intersection of cycles is well defined in intersection homology.
More precisely, if $\p$ and $\q$ are complementary perversities (this
means $\p + \q = \t$), there is a non degenerated bilinear map~:
$$IH^{\p}_i(X;{\bf Q})\times IH^{\q}_{n-i}(X;{\bf Q})\rightarrow
IH^{\t}_0(X;{\bf Q})
\buildrel \varepsilon \over \rightarrow {\bf Q}$$
corresponding to the intersection of cycles, followed by the evaluation
map $\varepsilon$.

If $\Sigma \subset X_{n-2}$, the Poincar\'e homomorphism, cap-product by
the fundamental class
$[X]$, admits the following factorisation, for every perversity $\p$~:
$$\matrix{
H^{n-i}(X) &&\hfleche\rightarrowfill{.\cap[X]}& &H_i(X)\cr
\vfleche\downarrow{\alpha} &&&&\vfleche\uparrow{\omega}\cr
IH^{\0}_i(X) &\rightarrow &IH^{\p}_i(X) &\rightarrow &IH^{\t}_i(X)
\cr}$$\par

Intersection homology theory is the good context to extend to singular
varieties results such that Morse theory [GM3], Lefschetz hyperplane
theorem [GM3],  hard Lefschetz theorem [BBD], Hodge decomposition [Sa]
and de Rham theorem (see the following paragraph).

\bigskip
\noindent{\bf 3. de Rham theorem for stratified varieties and polar
forms.} \medskip

The constructions we will use are taken from those of  Cheeger [Ch]
and Cheeger-Goresky-MacPherson [CGM], who proved  in particular
situations the ``standard" result concerning ${\cal L}^2$-cohomology of
differential forms,  i.e. the isomorphism~:
$$H^*_{(2)} (X-\Sigma) \cong {\rm Hom} (IH^{\p}_*(X,\RR );{\bf R})\ .$$
Many authors proved de Rham theorems for ${\cal L}^2$-forms, or ${\cal
L}^p$-forms, in different situations but always in the framework of
intersection homology (see [Bra] for a partial survey).

The  constructions that we give are also related to the theory of shadow
forms [BGM] which is another way to extend the de Rham theorem. In a
polyedron $(K)$ in the euclidean space ${\bf R}^n$, with a given
barycentric subdivision $(K')$, we associate, to each simplex $\sigma$
in $(K')$, a differential form $\omega(\sigma)$ in the interior of
simplices of maximal dimension, in a very explicit way. The shadow forms
have poles over faces of $(K)$~: If the defect of transversality of
$\sigma$ with a face $F$ of a $(K)$-simplex is $q$, then the maximum
order of poles of $\omega(\sigma)$ on this face is $q$. It can be proved
the inclusion quasi isomorphism~: $$\{ \hbox{ Shadow forms } \} \subset
\{\  {\cal L}^p - \hbox{ forms } \}$$

\bigskip

\goodbreak
\noindent{\bf The cone situation.}

We shall begin by defining an intersection complex of differential forms
on the cone $cL = [0,1[ \times L / \{ 0\} \times L$ with smooth
$(n-1)$-dimensional basis $L$ and vertex $\{ s\}$. This complex will
depend on two positive numbers the ``pinching number" $\alpha $ and the
``control number" $\beta $ of which we give now interpretations.

Recall that we want to characterize the cone by the behavior of
differentiable functions on its open regular stratum $]0,1[\times L$.
 The metric $dr^2+r^{2\alpha}  g_L$, where $g_L$ is a metric on the  link
$L$, separates cylinder ($\alpha =0$) and cones ($\alpha >0$), [Ch]. But
a metric is applied only to $k$-forms, $k>0$. To define a suitable action
at the functions level it seems natural to use, ``near the vertex" $\{
s\}$, the ``germ" action of the multiplicative group $\Rm$ on $cL$ given
by $\rho (r,x)=(\rho r,x)$ where $\rho \in \Rm$ and $(r,x)\in cL$.
Although this is not an isometry in the cone case (i.e. $\alpha >0$),
each $g\in C^\infty(L)$ determines a $1$-form $\omega =r^\alpha dg$
whose norm is equivariant~:
$\Vert\rho^*\omega\Vert =\Vert\omega \Vert$ and the functions  $r^\alpha
g$ are equivariant in $C^\infty(]0,1[\times L)$ relatively to the
modified germ action
$$(\rho f)(r,x)=\rho ^\alpha f({r\over\rho} ,x).$$
This action (and the associated equivariant relation) plays the role of a
metric at the functions level. But it does not respect the algebra
operation ($\rho ^\alpha$ acts as a derivation, cf. (1) of introduction).
To solve this problem we remark that the $r$-derivatives of the functions
$r^\alpha g$ verify the equivariant condition~:
$$\hbox{for any } n\in {\bf N},\quad \vert r^{-\alpha +n }
(r\partial_r)^n(r^\alpha g)\vert \quad \hbox{ is independant of } r\ .$$
So we modify the $C^\infty$-topology of
$C^\infty(]0,1[\times L)\cong C^\infty(]0,1[)\widehat\otimes C^\infty(L)$
using the following semi-norms on $C^\infty(]0,1[)$
$$\sup_{r\in ]0,1[}\vert r^{-\alpha +n }(r\partial_r)^n( \ - \ )\vert$$
This motives the introduction of the ``reference algebra"
that is the Frechet algebra of the (near $\{s \}$) bounded smooth
functions relatively to the multiplicative group $\Rm$~:
$$A=A(r)=\{a\in
C^\infty(]0,1[)  \ :\  \forall \ k \ \sup_{r\in ]0,1[}\vert
(r\partial_r)^k a(r)\vert<+\infty\}.$$
Then we substitute equivariant functions by controlled functions, which
are boun\-ded relatively to these semi-norms, i.e.
the elements of
$$r^\alpha  A\widehat\otimes C^\infty(L)$$
where $\widehat\otimes $ is the completed projective tensor product.
The lack of algebra structure coming from the equivariant framework
disappears when considered in the Frechet algebra framework.
\par

The number $\beta $ controls the order of  poles of the forms near the
vertex $\{ s\}$.
The differential $A$-module of poles is~:
$$M^*_\beta = r^{-\beta} A \oplus  r^{-\beta} A {dr\over r}$$
Remark that for $f_i \in C^\infty(L)$, the equivariant k-form
$\omega _k=r^{-\beta }(r^\alpha f_0)(r^\alpha df_1) \wedge \cdots \wedge
(r^\alpha df_k)$ has a pole of order $\beta -(k+1)\alpha$ in $\{ s\}$.
\medskip

Now we shall construct the intersection complex. The $A$-module of
$\beta$-controlled differential forms on the cone $cL$  is defined by~:
$${\cal B} _\beta ^k(cL)= r^{(k+1)\alpha} \left [M^*_\beta
\widehat\otimes \Omega^*(L)\right ]^k$$
Let $\{ U_i \}_{i\in I}$ be a locally finite atlas of $L$, we
denote by $x=\{x_1,\ldots,x_{n-1} \}$ a system of local coordinates in
$U_i$.
For $\gamma\in {\bf R}$, it can be shown that the $A$-module
$r^\gamma A\widehat\otimes C^\infty(L)$ of $\gamma$-controlled functions
on $cL$, denoted by $C^\infty_{\gamma}(cL)$, verifies~:
$$\eqalign{
C^\infty_{\gamma}(cL)=\{f
\in {\cal C}^\infty(]0,1[ \times L)\ :&\ \forall U_i\ , \forall s \ ,
\forall \lambda=(\lambda_1,\ldots,\lambda_{n-1}), \cr
&\sup_{(r,x) \in ]0,1[ \times U_i} r^{\gamma + s} \Bigl\vert
{\partial^{s+\vert
\lambda \vert}f(r,x) \over \partial r^s \partial^\lambda x}\Bigr\vert
<\infty\}\cr}$$
where $s\in {\bf N}$, $\lambda_i \in  {\bf N}$,
$\partial^\lambda x =
(\partial x _1)^{\lambda_1} \cdots (\partial x_{n-1})^{\lambda_{n-1}} $
and $\vert\lambda \vert =\lambda_1+\ldots+\lambda_{n-1}$.

Then ${\cal B} _\beta ^k(cL)$ is the module of differential forms
$\omega \in \Omega^k(]0,1[\times L)$ whose restriction, in each $U_i$, is
a sum of elements of the type $ a dx^k + b{dr\over r}\wedge dx^{k-1} $
with $a,b\in C^\infty_{(k+1)\alpha-\beta} (cL)$. It is not a complex (cf.
(2) in introduction and also \S 2) so we define the
intersection complex by~:  $$I{\cal B} _\beta ^k(cL)=\{\omega
\in {\cal B} _\beta ^k(cL) \ : \ d\omega \in {\cal B} _\beta
^{k+1}(cL)\}$$ With minor changes of control parameters, the following
result  is similar to [BL], Th\'eor\`eme 2.4.

\noindent {\bf Theorem 0}. {\sl Let $cL$ be a cone over an
$(n-1)$-dimensional smooth manifold $L$, and $\bar p$ any perversity
such that $p_n = n-2 - [\beta /\alpha ]$. Suppose that $\beta/\alpha$ is
not an integer, then}
$$H^k(I{\cal B} _\beta ^*(cL)) \cong
\hbox{Hom}(IH_k^{\bar p}(cL,{\RR});{\bf R})=
\cases{H^k_{\rm dR}(L)&if $k\le[\beta /\alpha ]-1$\cr
0&if $k>[\beta /\alpha ]-1$\cr}$$

\goodbreak
The idea of the proof  is the following (cf [Ch])~: firstly by
Poincar\'e Lemma, pointed in the vertex $\{ s\}$, all closed forms on
$cL-\{s\}$ are cohomologous to the extension (constant in $r$) to
$cL-\{s\}$ of a closed form on $L$. Now~:\par
\noindent - on one hand, such an extension is controlled only if $k\le
[\beta /\alpha ]-1$, \par
\noindent - on the other hand, for high degrees the order of poles of
controlled forms discreases with the degree of the form. Namely, for $k>
[\beta /\alpha ]-1$, the form converges to 0 as $r$ goes to 0.
\bigskip

\noindent {\bf Atlases of iterated cones and $\bar \beta $-controlled
forms.}

In order to generalize the result to stratified varieties, we will use
atlases whose charts are iterated cones.

Let $x$ be a point in a stratum $S_{n-j_1}$ and $\psi_x:U_x\buildrel
\cong\over\to  B^{n-j_1}\times cL_x$ a distinguished open neighborhood
as previously described. The link $L_x$ is a singular variety and is
covered by
distinguished open sets of the same type. By iteration, we obtain a
chart which defines an {\sl iterated cone}~:
$$W_{\jbar} =  {\bf B}^{n-j_1}\times c\left({\bf B}^{j_1-1-j_2} \times
c\left({\bf B}^{j_2-1-j_3} \times \cdots \times c({\bf B}^{j_{\ell}-1})
\cdots
\right) \right)$$
where ${\jbar} = \{n+1= j_0>j_1>j_2>\cdots>j_{\ell}>j_{\ell +1}=0\}$
denotes a decreasing sequence of integers and ${\bf B}^{j_t-1-j_{t+1}}$
is an open ball in ${\bf R}^{j_t-1-j_{t+1}}$ (possibly a point).

Via the homeomorphism $\psi_x$,  this chart corresponds to the following
chain of elements of the filtration of $X$~:
$$\emptyset=X_{n-j_0}\subset X_{n-j_1}
\subset X_{n-j_2}\subset
\cdots \subset X_{n-{j_\ell}} \subset X=X_n=X_{n-j_{\ell+1}}\ .$$
We will denote in the same way the iterated cone and its image in $X$.

We obtain, by this way, a covering of $U_x$ by charts which are iterated
cones and corresponding to different sequences $\jbar$.

Let us describe the coordinates in $W_{\jbar}$. For $t=0,\ldots,\ell$,
we denote by $u^t_{a_t}$ the coordinates in ${\bf B}^{j_t-1-j_{t+1}}$
($1\le a_t \le j_t-1-j_{t+1}$) and $r_{j_t}$ the coordinate of the
generatrix of the $t$-th cone. So, we have $\ell$  coordinates of the
type $r_{j_t}$ and $n-\ell$ coordinates of the  type $u^t_{a_t}$. We
set~:

$\bar r = (r_{j_1},\ldots ,r_{j_\ell})$ \quad
$\bar u =({\tilde u}_{j_0},\ldots,{\tilde u}_{j_\ell})$,
with ${\tilde u}_{j_t}=(u^t_{1},\ldots,u^t_{j_t-1-j_{t+1}})$

Given such an atlas on a Thom-Mather stratified space $X$, we can define
a metric, called $\bar\alpha$-metric on $X-\Sigma$, and more precisely on
iterated  cones, in the following way~: consider a sequence of real
numbers $\bar\alpha = (\alpha_0,\cdots,\alpha_n)$ associated to the
filtration $(*)$, each $\alpha_j$ corresponding to the stratum $S_{n-j}$.
On each open ball $B^{k}$, with coordinates $u_1,\ldots,u_{k}$, the
metric is the euclidean one~: $du_1^2 + \cdots + du_{k}^2$. On the
regular part of each product  $B^{j_t-1-j_{t+1}}\times c(L)$, the metric
is the product metric $(du^t_1)^2  + \cdots + (du^t_{j_t-1-j_{t+1}})^2
+ (dr_{j_{t+1}})^2 + (r_{j_{t+1}})^{2\alpha_{j_{t+1}}}g_L $
where $r_{j_{t+1}}$ is the coordinate of the generatrix of the
cone and $g_L$ is the metric on the regular part of $L$, defined
inductively.

Now we can define controlled functions on iterated cones. Let $\bar
\gamma =(\gamma _1,\cdots,\gamma _n)$ an $n$-uple of real numbers, we
define $$
C^\infty_{\bar \gamma}(W_{\jbar}) =
\left( \widehat{\otimes}_t\  r^{\gamma_{j_t}}_{j_t} A(r_{j_t}) \right)
\widehat\otimes
\left( \widehat{\otimes}_t\  C^{\infty}({\bf B}^{j_t-1-j_{t+1}}) \right)
$$
In an equivalent way, this is the set of functions $f\in C^\infty
(W_{\jbar} \cap (X-\Sigma))$ such that for all $n$-uples of positive
integers  $(\bar s,\bar\lambda)=(s_1,\ldots,s_\ell,\bar
\lambda_0,\ldots,\bar\lambda_{\ell})$,
with $\bar\lambda_t = ({\lambda^t_1},\ldots,{\lambda^t_{j_t-1-j_{t+1}}})
$, $$
\sup_{(\bar r ,\bar u )} (\bar r)^{\bar\gamma+\bar s} \left\vert{
\partial^{\vert (\bar s,\bar\lambda) \vert}f(\bar r,\bar u) \over
(\partial\bar r)^{\bar s} (\partial \bar u)^{\bar\lambda }}\right\vert
< + \infty
$$
where $(\bar r)^{\bar\gamma}=(r_{j_1})^{\gamma_{j_1}}\cdots
(r_{j_\ell})^{\gamma_{j_\ell}}$ and  $(\partial\bar r)^{\bar s} =
\partial(r_{j_1})^{s_1} \cdots \partial(r_{j_\ell})^{s_\ell}$, and in
the same way $(\partial \bar u)^{\bar\lambda }= (\partial \tilde
u_{j_0})^{\bar{\lambda_0}}
\cdots (\partial \tilde u_{j_\ell})^{\bar{\lambda_\ell}}$ where
$(\partial \tilde u_{j_t})^{\bar{\lambda_t}} = (\partial
u^t_{1})^{\lambda^t_1}
\cdots (\partial u^t_{j_t-1-j_{t+1}})^{\lambda^t_{j_t-1-j_{t+1}}}$.
\medskip

Let $U$ be an open set in $X$, we say that a function $f\in C^\infty
(U-\Sigma) $ is  $\bar\gamma$-controlled on $U$ and we denote $f\in
C^\infty_{\bar \gamma}(U)$  if, for all $W_{\jbar}$, we have $f\in
C^\infty_{\bar \gamma}(W_{\jbar}\cap U)$. The cor\-res\-pon\-dance $U
\mapsto C^\infty_{\bar \gamma}(U)$ defines a presheaf, which is not a
sheaf. The associated sheaf, independant of the atlas and denoted by
${\cal C}^\infty_{\bar \gamma}$ is defined by
$${\cal C}^\infty_{\bar \gamma}(U)=\{ f\in C^\infty(U-\Sigma )  \  :
\forall \ x\in U-\Sigma ,
\exists  \ V_x\subset U-\Sigma , f\in C^\infty_{\bar \gamma}(V_x)\}.$$

Remark that for each open set $U$ such that $\{ s\} \in U$, then
$C^\infty_{\bar \gamma}(U)=\Gamma(U,{\cal C}^\infty_{\bar \gamma})$.

A differential form  defined in $W_{\jbar}\cap U-\Sigma$ can be
written as a sum of terms of the form~:
$$a(\bar r, \bar u) (d\bar r)^{\bar \mu}\wedge (d\bar u)^{\bar k}$$
where $(d\bar r)^{\bar \mu}= (dr_{j_1})^{\mu_1} \wedge \cdots \wedge
(dr_{j_\ell})^{\mu_\ell}$, $(\mu_i = 0,1)$,
$(d\bar u)^{\bar k}= (d{\tilde u}_{j_0})^{\bar{k_0}}\wedge \cdots
\wedge (d{\tilde u}_{j_\ell})^{\bar{k_\ell}}$ and $(d{\tilde u}_{j_t})
^{\bar{k_t}} = (du_{1}^{t})^{k^t_{1}}\wedge \cdots \wedge
(du_{j_t-1-j_{t+1}}^{t})^{k^t_{j_t-1-j_{t+1}}}$. Here, $k_t$ is the
number of coordinates appearing in $d{\bar u}_{j_t}$, i.e. $k_t=\vert
\bar{k_t}\vert=k^t_1 + \cdots k^t_{j_t-1-j_{t+1}}$.\par

Let $\bar\beta = (\beta_1,\ldots,\beta_n)$ be a fixed sequence of
strictly positive real numbers.

Let $\omega$ be a differential form defined on $U-\Sigma$, we say that
$\omega$ is  ${\bar\beta}$-controlled on $U$, relatively to the
$\bar\alpha$-metric, if for all $W_{\jbar}$ all the coefficients $a(\bar
r,\bar u)$ belong to $C^\infty_{\bar \gamma}(W_{\jbar} \cap U)$ where,
for all $1\le t\le \ell$, $\gamma_{j_t}= \beta_{j_t}+\mu_t +
(k_t+1)\alpha_{j_t}$. In this expression, $\alpha_{j_t}$ is the pinching
of the corresponding cone, $\beta_{j_t}$ determines the order of pole,
and $\mu_t$ and $k_t$ are determined by $\omega$.\goodbreak

Let us denote by $I{\cal B}_{\bar \beta}^k(U)$ the
space of $k$-differential forms $\omega$ such that $\omega$ and
$d\omega$ are ${\bar\beta}$-controlled on $U$.
We define the sheaf complex of ${\bar\beta}$-controlled differential
forms ${\cal IB} _{\bar \beta}^*$ as the sheaf associated to the
presheaf $U\mapsto I{\cal B} _{\bar \beta}^k(U)$ in the same way that
we defined ${\cal C}^\infty_{\bar \gamma}$.  \par

The proof of the following theorem  is similar to [BL],   Th\'eor\`eme
3.5, up to minor changes of control parameters.
\medskip
\noindent {\bf Theorem 1.}
{\sl Let $X$ be a Thom-Mather stratified space, endowed with a covering
by iterated cones and an $\bar\alpha$-metric. Suppose  $\bar \beta$
given, satisfying the following perversity condition~:
$$\left[{\beta_j\over \alpha_j}\right] \leq \left[{\beta_{j+1}\over
\alpha_{j+1}}\right] \leq\left[{\beta_j\over \alpha_j}\right] + 1
\qquad \forall j\ ,$$
with $\displaystyle{\beta_j/\alpha_j}$ non
integer. Then, there is an isomorphism
$$H^k(I{\cal B} _{\bar \beta}^*(X)) \cong
\hbox{Hom}(IH_k^{\bar p}(X,{\RR});{\bf R})$$
with $\bar p_j=j-2-{\displaystyle \left[\beta_j \over \alpha_j\right]}$}

\medskip
Remark~: The following inclusions are quasi-isomorphisms (the first one
being defined only for polyedra)~:
$$\hbox{ $\{$ Shadow forms $\} \subset \{ \ {\bar\beta}$-controlled forms
$\} \subset \{ {\cal L}^p$- forms$\}$}$$
\bigskip
\goodbreak
\noindent{\bf 4. Hochschild and cyclic homology of controlled functions.}
\medskip

For the rest of the paper and for simplicity, we turn back to the case
of a cone over a smooth manifold. Firstly we show that controlled
functions generate the intersection complex (Theorem 2), then we give
relation with Hochschild and cyclic homology (Theorem 3).\par

\noindent{\bf How to generate controlled forms using controlled
functions.}

We remark that there are two methods to obtain a complex, starting
with ${\cal B}_\beta^*$ (cf. (2), introduction)~: in the first one,  we
consider the intersection complex $I{\cal B}_\beta^*(cL)$ previously
defined, in the second one we stabilize by the de Rham operator, i.e. we
add coboundaries. It is not difficult to show directly that the two
complexes are quasi-isomorphic.  The second one has the important
property to be generated by controlled functions. \par

Firstly, dealing with a complex (namely
$\omega$ and $d\omega$ are ${\bar\beta}$-controlled) is translated in
cyclic homology theory by the use of unitarized algebras. We define the
$A$-unitarization of the $A$-algebra $ r^{\alpha} A \widehat \otimes
C^\infty(L)$ as the algebraic sum in  $C^\infty(]0,1[\times L)$
$$IC^\infty_\alpha(cL) = r^{\alpha} A \widehat \otimes C^\infty(L) + A $$
The Frechet $A$-algebra structure is provided when we identify
$IC^\infty_\alpha(cL)$ with the quotient algebra $( r^{\alpha} A \widehat
\otimes C^\infty(L)\oplus A)/I$ where $I=\{f+a \ : \forall \ g \in
r^{\alpha} A \widehat \otimes C^\infty(L), (f+a)g=0\}$ is a closed ideal.
\par

Every differential form $\omega \in \Omega^*(]0,1[\times L)$ can be
written in an unique way  $\omega=\eta + {dr\over r} \wedge\varphi $
where $\varphi = i_{\scriptstyle{\partial\over\partial_r}} \omega $ and
$\eta\in \Omega ^*_P$, the space of differential forms relatively to the
projection $P:]0,1[\times L\to ]0,1[$.\par

Let us denote by $\Omega^*(IC^\infty_\alpha(cL))$ the subcomplex of
$\Omega^*(]0,1[\times L)$ generated by the functions which belong to
$IC^\infty_\alpha(cL)$. With the previous notation,
$\Omega^*_P(IC^\infty_\alpha(cL))$ is the complex  generated by the
element $\eta$ in $\Omega^*(IC^\infty_\alpha(cL))$ and with differential
$d_L$ induced by the de Rham differential $d=d_L+d_r$. If the cone is
looked as a family of spaces, $L_r=L$ for $r\in]0,1[$ and $L_0=\{s\}$,
(this point of vue appears implicitely in the $A$-module structure),
then $\Omega^*_P(IC^\infty_\alpha(cL))$ is the complex of sections of
the family of de Rham complexes on $]0,1[$ associated to the $A$-algebra
$IC^\infty_\alpha (cL)$.\par

Using the previous definition of the $A$-module of poles $M^*_\beta $ we
define the complex
$$I\Omega^*_\beta(cL) =M^*_\beta\widehat\otimes_A
\Omega^*_P(IC^\infty_\alpha(cL))$$
with differential
$$\Diag{r^{-\beta} \ \Omega_P^k(IC^\infty_\alpha(cL))&\oplus & r^{-\beta}
{\displaystyle{dr\over r}}\wedge \  \Omega_P^{k-1}(IC^\infty_\alpha(cL))
\cr \downarrow d_L&\searrow d_r&\downarrow d_L\cr
r^{-\beta} \ \Omega_P^{k+1}(IC^\infty_\alpha(cL))&\oplus & r^{-\beta}
{\displaystyle{dr\over r}}\wedge \
\Omega_P\Omega^k(IC^\infty_\alpha(cL))\cr }$$
We observe that the problem (2) of introduction is solved and that the
unitarization with ${\bf R}$ would not be sufficient to obtain a complex.

\noindent {\bf Theorem 2.} {\sl There is an isomorphism of complexes~:
$$I\Omega^*_\beta(cL)   \cong {\cal B}_\beta^* + d{\cal B}_\beta^*$$
In the assumptions of the Theorem 0, we have
$$H^k(I\Omega {\cal B} _\beta ^*(cL))=
\hbox{Hom}(IH_k^{\bar p}(cL,{\RR});{\bf R}).$$}
\noindent{\bf Hochschild and periodic cyclic homology of controlled
functions.}

Now, let us give the relation with Hochschild and cyclic homology.
Firstly we recall some basic and general properties, the references are
[Co] and [Lo]. Let $\Lambda$ be a field and ${\cal A}$ be an algebra
with unit, so $\Lambda \subset{\cal A}$. The Hochschild complex
$(C_*({\cal A}),b)$ is defined by $C_k({\cal A})={\cal A}\otimes {\cal
A}^{\otimes^k}$ and the Hochschild boundary is
$$b(a_0\otimes \cdot\cdot\otimes a_k) = \sum_{j=0}^{k-1}
(-1)^j a_0\otimes\cdot\cdot\otimes a_ja_{j+1}\otimes\cdot\cdot \otimes
a_k \ + \ (-1)^k a_ka_0\otimes a_1\otimes\cdot\cdot\otimes a_{k-1}$$
Its homology, called Hochschild homology, is denoted by $HH_*({\cal A})$.
The reduced Hoch\-schild complex $(C_*^{red}({\cal A}),b)$  is the
quotient of the Hochschild complex by the subcomplex generated by the
elements $a_0\otimes \cdots\otimes a_k$ where $a_i\in \Lambda$ for some
$i>0$. The reduced Hoch\-schild complex is quasi-isomorphic to the
Hochschild complex.

The Hochschild complex is a cyclic module, i. e. it admits a cyclic
action $$\tau (a_0\otimes \cdots\otimes a_k)=(-1)^k
a_k\otimes a_0\otimes a_1\otimes\cdots\otimes a_{k-1}\ .$$
We have $\tau ^{k+1}=id$, so $\tau $ defines an action of
${\bf Z}/(k+1){\bf Z}$.

The Connes cyclic homology is defined in the following way~: Consider
the situation where $\Lambda$ is a field of characteristic $0$. The
cyclic homology of ${\cal A}$, denoted by $HC_*({\cal A})$, is the
homology of the complex $(C_*({\cal A})/(1-\tau ),b)$ where $b$ is
induced by the Hochschild boundary. The relation between Hochschild and
Connes homology is given by the Connes exact sequence
$$ \cdots\to HH_k({\cal A})\buildrel I\over{\to }HC_k({\cal A})\buildrel
S\over{\to }
HC_{k-2}({\cal A})\buildrel B\over{\to }HH_{k-1}({\cal A})\to\cdots$$
Using the so called {\sl periodicity operator} $S$ we define the
periodic cyclic homology
$$PHC_*=\lim_k \left [ HC_k({\cal A})\buildrel S\over{\to } HC_{k-2}
({\cal A}) \right ].$$
The importance of the above definitions appears with the following
result. Let $X$  be a compact ${\cal C}^\infty$ manifold, ${\cal
A}=C^\infty(X)$ the Frechet algebra of differentiable functions on $X$
and $\Omega ^*(X)$ the associated de Rham algebra. Replace everywhere
$\otimes $ by the projective tensor product $\widehat\otimes$. Then
theapplication $\pi : C_k(C^\infty(X))\to \Omega ^k(X)$ defined by
\hfill\break
$\pi(f_0\otimes\cdots\otimes f_k)=f_0df_1\wedge\cdots\wedge df_k$
induces the following isomorphims [Co]~:
$$HH_*(C^\infty(X))\cong \Omega^*(X) \ ; \ PHC_{{\rm odd}\atop{\rm even}}
(C^\infty(X))\cong
\oplus_{{\rm odd}\atop{\rm even}}H^*_{dR}(X)$$
\hskip 23pt If $\Lambda$ is a
ring, we need a more general setting, the notion of {\sl mixed complex}
that we briefly  describe, as it appears in the singular framework.

A mixed complex, [Ka], $(M_*,b,B)$ is a graded module with two
differentials, $b$ of degree $-1$ and $B$ of degree $+1$ such that
$bB+Bb=0$. It defines a bicomplex $M_*[u]$ with differentials
$b(m u^k)=(bm)u^k$, $B(mu^k)=(Bm)u^{k-1}$ where degree$(u)=2$.
We define the {\sl Hochschild homology} of $(M_*,b,B)$ as $H_*(M_*,b)$
and the {\sl cyclic homology} as $H_*(M_*[u],b+B)$. There is again a
Connes exact sequence and we can define the periodic cyclic homology.

When $\Lambda$ is a field of characteristic $0$, the relation between the
two previous definitions is the following. Replacing the quotient
$C_k({\cal A})/(1-\tau )$ by a ${\bf Z}/(k+1){\bf Z}$-free resolution of
$C_k({\cal A})$, we obtain a bicomplex which is quasi-isomorphic to the
bicomplex associated to a mixed complex. The mixed complex structure of
$C_*({\cal A})$ is given by the  Hochschild boundary $b$ and by the
operator  $B$ defined by
$$\eqalign{B(a_0\otimes\cdots\otimes a_k)=&
\sum_{j=0}^{k-1}(-1)^{kj}1\otimes a_j\otimes\cdots\otimes a_k\otimes
a_0\otimes\cdots \otimes a_{j-1}\cr
&-(-1)^{k(j-1)}a_{j-1}\otimes 1\otimes a_j\otimes \cdots\otimes a_k
\otimes a_0\otimes\cdots \otimes a_{j-2}\ .\cr}$$
Then the two definitions of cyclic homology agree.
\goodbreak
\medskip

In the following, the cone must be seen as a family of spaces $L_r=L$
for $r>0$ and $L_0=\{ s\}$ and we shall use a slight generalization of
the previous cyclic construction.

Consider  the  $A$-Hochschild complex
$$C_k^A(IC^\infty_\alpha(cL))=IC^\infty_\alpha(cL) \widehat\otimes_A
\cdots\widehat \otimes_A IC^\infty_\alpha(cL)$$
($(k+1)$-terms)  and denote by $b_A$ its differential.
Define the Hochschild-intersection complex by
$$IC_k^\beta(cL)        =M^*_\beta\widehat\otimes_A
C_{k-*}^A(IC^\infty_\alpha(cL))$$
The elements of degree $k$ are  sum  of terms
$$r^{-\beta }f_0 \otimes f_1\otimes \cdots\otimes f_k
+ r^{-\beta }{dr\over r}g_0\otimes g_1\otimes
\cdots\otimes g_{k-1}$$
such that $f_i, g_j\in   IC^\infty_\alpha(cL)$. The total differential
is given by
$$\Diag{r^{-\beta} \  C_k^A(IC^\infty_\alpha(cL))&\oplus &
r^{-\beta} A{\displaystyle{dr\over r}}\widehat\otimes_A\
C_{k-1}^A(IC^\infty_\alpha(cL))\cr
b^{k+1}_A\uparrow\downarrow B^k_A
&\searrow d_r&b^k_A\uparrow\downarrow B^{k-1}_A\cr
r^{-\beta} \  C_{k+1}^A(IC^\infty_\alpha(cL))&\oplus &
r^{-\beta} A{\displaystyle{dr\over r}}\widehat\otimes_A \
C_k^A(IC^\infty_\alpha(cL))\cr }$$
where the operator $B_A$ is associated to the cyclic operation (as in
the classical  case, [Lo])
$$\tau(r^{-\beta } g_0\otimes\cdots\otimes g_k)= (-1)^k
r^{-\beta }g_k\otimes g_0\otimes\cdots\otimes g_{k-1}$$
(we leave the factor $r^{-\beta }$ in the first term cf. (3),
introduction) and $d_r$ corresponds to the r-derivation,
$$\eqalign{d_r(r^{-\beta }f_0\otimes f_1\otimes \cdots\otimes f_k)=\cr
(-1)^kr^{-\beta }{dr\over r}&\wedge\left[-\beta f_0\otimes\cdots\otimes
f_k +\sum_{i=0}^k f_0\otimes\cdots\otimes \partial_r f_i\otimes
\cdots\otimes f_k\right]\cr}$$
So $b^*_A\oplus b^{*-1}_A$ has degree $-1$ and
$B^*_A\oplus B^{*-1}_A  \ + \ d_r$ has degree $1$.
\medskip

\noindent {\bf Lemma.} {\sl The triple
$$(IC_*^\beta(cL),b=b^*_A\oplus b^{*-1}_A, B=B^*_A\oplus B^{*-1}_A  \ +
\ d_r)$$
is a mixed complex
(i. e. $b^2=B^2=bB+Bb=0$).}
\medskip

\noindent {\bf Theorem 3.} {\sl i) The Hochschild homology of
$IC_*^\beta(cL)$
(with differential $1\otimes b^*_A$)  is~:
$$HH_k(IC_*^\beta(cL))  \cong I\Omega^k_\beta(cL)$$
ii) the periodic cyclic homology is~:
$$PHC_k(IC_*^\beta(cL)) \cong IH^k_{\bar q}(cL)\oplus
IH^{k-2}_{\bar q}(cL) \oplus \cdots$$
where the perversity $\bar q$ satisfies $ q_n=[{\beta\over \alpha}]-1$}

\noindent {Sketch of the proof}~:
i) The terms in $r^{-\beta}$ do not modify the demonstration (they stay
as common factor), so we omit  them in the proof of part (i). Consider
the  reduced
$A$-Hochschild complex
$$\eqalign{C_k^{\rm red}&=IC^\infty_\alpha(cL)\widehat \otimes_A
 C^\infty_\alpha(cL)
\widehat\otimes_A  \cdots \widehat\otimes_A C^\infty_\alpha(cL)\cr
&\cong r^{k\alpha }IC_\alpha^{\infty}(c(L))\widehat\otimes
C^\infty(L^{\times k})
\subset C^\infty(c(L^{\times (k+1)})-\{s\})\ .\cr}$$
Using the lemma below, it suffices to prove that $H_k(C_*^{\rm red})
\cong \Omega_P^k(IC^\infty_\alpha(cL))$.

For every open $U\subset [0,1[$ we can define the Frechet module of
controlled functions on $U\times L^{\times (k+1)}$ in the same way, as
$C_k^{\rm red}$.  So we have  a presheaf $U \mapsto {\cal C}_k^{\rm red}
(U)$ which define a fine sheaf ${\cal C}_k^{\rm red}$ using the same
localization condition as ${\cal C}_{\bar \gamma }^\infty$. Its space of
sections is the reduced $A$-Hochschild complex $C_k^{\rm red}= {\cal
C}_k^{\rm red}([0,1[)$. We can also associate to
$\Omega_P^k(IC^\infty_\alpha(cL))$  a fine sheaf which is denoted by
${\cal I}\Omega^k $. We define a sheaf morphism $\pi :  {\cal C}_k^{\rm
red}\to {\cal I}\Omega^k $ as above~: for each $U$ set $\pi (f_0\otimes
f_1\otimes\cdots\otimes f_k)=f_0df_1\wedge\cdots\wedge df_k$ where $f_i$
is controlled on $U$. For each $r\in ]0,1[$, the fiber complex $({\cal
C}_*^{\rm red})_r$ is isomorphic to the standard Hochschild complex of
$C_*(C^\infty(L))$ and the fiber complex $({\cal I}\Omega^*)_r $ is
isomorphic to $\Omega^*(L)$. These isomorphisms are induced by the
following one in degree 0~: $$\eqalign{IC^\infty_\alpha(cL)_r &\ \
\longrightarrow\quad C^\infty(L)\cr r^\alpha a \otimes f + c
&\quad\mapsto\quad r^\alpha a  f + c  \cr}$$ By Connes's theorem [Co],
for $r>0$,  the Hochschild homology of the fiber complex of ${\cal
C}_*^{\rm red}$ is isomorphic to $\Omega^*(L)$ so the morphism $\pi $
induces a fiber isomorphism between the homology sheaf  ${\cal H}_k(
{\cal C}_*^{\rm red})$ and the sheaf ${\cal I}\Omega^* $. This can be
extended for $r=0$ and this implies that ([Go], 4.5)~: $$H_k(\Gamma
(cL,{\cal C}_*^{\rm red})) \cong \Gamma(cL,{\cal H}_k({\cal C}_*^{\rm
red})) \ .$$

ii) By the isomorphism $H_k(C_*^{\rm red}) \cong \Omega_P^k
(IC^\infty_\alpha(cL))$, the differential $B$ gives the de Rham
differential. Using theorems 1 and 2, the proof of ii) is then similar
to  the non singular case. \medskip

Similar sheaf arguments give the proof of the following lemma (used in
the previous demonstration)~: \medskip
\noindent {\bf Lemma}. {\sl The $A$-Hochschild complex
$C_*^A(IC^\infty_\alpha(cL))$
and the reduced $A$-Hochschild complex $C_*^{\rm red}$  are
quasi-isomorphic.}
\bigskip
\goodbreak
\noindent {\bf Bibliography}\moyen

\parindent = 0pt

\item {[BBD]} A. Beilinson, J. Bernstein et P. Deligne. {\sevensl
Faisceaux pervers}, Asterisque n${}^{\rm o}$ 100, 1983.

\item {[BG]} J. Block and E. Getzler. {\sevensl Equivariant cyclic
homology and equivariant differential forms}, Ann. Scient. \'Ec. Norm.
Sup. t. 27, (1994), 493 - 527.

\item {[Bra]} J.P. Brasselet. {\sevensl De Rham's theorems for singular
varieties},  Contemporary Mathematics 161, (1994), 95 - 112.

\item {[BGM]} J.P. Brasselet, M. Goresky and R. MacPherson. {\sevensl
Simplicial Differential Forms with Poles},  Amer. Journal of Maths.,
113 (1991), 1019-1052.

\item {[BL]} J..P. Brasselet et A. Legrand. {\sevensl Un complexe de
formes diff\'erentielles \`a croissance born\'ee sur une vari\'et\'e
stratifi\'ee}, Annali della Scuola Normale Superiore di Pisa Serie IV,
Vol XXI, Fasc.2 (1994), p. 213 - 234

\item {[Bry]} J.L. Brylinski. {\sevensl Cyclic homology and equivariant
theories}, Ann. Inst. Fourier, vol. 37, 4 (1987), 15 - 28.

\item {[Ch]} J. Cheeger. {\sevensl On the Hodge theory of Riemannian
pseudomanifolds,} Proc. of Symp. in Pure Math. vol 36 (1980), 91-146.
Amer.Math.Soc., Providence R.I.

\item {[CGM]} J. Cheeger, M. Goresky and R. MacPherson.
{\sevensl ${\cal L}^2$-cohomology and intersection cohomology for
singular varieties,} Seminar on Differential Geometry, S.T.\-Yau, ed.
Ann. of Math. Studies, Princeton  University Press, Princeton N.J., 102
(1982), 303-340.

\item {[Co]} A. Connes. {\sevensl Non commutative differential
geometry,} Publ. Math. I.H.E.S., 62 (1986),  257 - 360.

\item {[Go]}  R. Godement. {\sevensl Th\'eorie des faisceaux},
Actualit\'es scientifiques et industrielles, 1252, Hermann 1964.

\item {[GM1]} M. Goresky and R. MacPherson. {\sevensl Intersection
homology theory}, Topology 19 (1980), 135--162.

\item {[GM2]} M. Goresky and R.MacPherson. {\sevensl Intersection
homology theory II,} Inv. Math. 71 (1983), 77--129.

\item {[GM3]} M. Goresky and R. MacPherson. {\sevensl Stratified Morse
theory}, Ergebnisse. Band 14, Springer-Verlag 1987.

\item {[Ka]} C. Kassel. {\sevensl Cyclic homology. Comodules and mixed
complexes,} J. of Algebra 107 (1987), 195--216.

\item {[Lo]} J.L. Loday. {\sevensl Cyclic Homology,} Grund. der math.
Wiss. 301, Springer Verlag (1992)

\item {[Ma]} J. Mather. {\sevensl Notes on topological stability},
Harvard University, 1970.

\item {[Sa]} M. Saito. {\sevensl Modules de Hodge polarisables}, IAS,
Princeton, 1986.

\bigskip
\bigskip
\hbox to \hsize{\hfill
{\moyen
\parskip 0pt
\baselineskip=12pt
\vtop {\vskip 0.5truecm} {\hskip 0.1truecm}
\vtop{\hsize=0.5\hsize
\obeylines{
Jean-Paul Brasselet
IML - CNRS
Luminy Case 930
F-13288 Marseille Cedex 9
e-mail : jpb@iml.univ-mrs.fr
}}\hfill
\vtop{\hsize=0.45\hsize
\obeylines{
Andr\'e Legrand
Laboratoire Emile Picard
Universit\'e Paul Sabatier
118 Route de Narbonne
F-31062 Toulouse Cedex
e-mail : legrand@picard.ups-tlse.fr
}}}\hfill}

\end